\documentclass{scrartcl}
\usepackage[english]{babel}
\usepackage[latin1]{inputenc}
\usepackage{amsthm,amsfonts,amssymb,amsmath}
\usepackage{graphicx,color}

\usepackage[pdftex]{hyperref}
\usepackage{textcomp}
\usepackage{multirow}
\usepackage{mathrsfs}
\usepackage[normalem]{ulem}
\usepackage{indentfirst}
\usepackage{geometry}

\begin{document}

\title{On the mistake in defining fractional derivative
using a non-singular kernel}
\author{E. Capelas de Oliveira\thanks{capelas@ime.unicamp.br}
\and 
S. Jarosz\thanks{stjarosz@gmail.com}
\and 
J. Vaz Jr.\thanks{vaz@unicamp.br}} 
\date{}
\publishers{Department of Applied Mathematics\\
University of Campinas\\
Campinas, SP, Brazil}

\maketitle

\begin{abstract}
Definitions of fractional derivative of order $\alpha$ ($0 < \alpha \leq 1$) 
using non-singular kernels have
been recently proposed. In this note we show that these
definitions cannot be
useful in modelling problems with 
an initial value condition (like, for example, 
the fractional diffusion equation) 
because the solutions obtained for 
these equations do not satisfy 
the initial condition (except for the integer case 
$\alpha = 1$). In order to satisfy 
an arbitrary initial condition the definitions of 
fractional derivative must 
necessarily involve a singular kernel.   
\end{abstract}

Fractional calculus has sparked the interest from researchers ever since its
beginning, but especially in recent decades, possibly because of its
many applications, like in the modelling of memory effects or anomalous diffusion. 
One of the key concepts of fractional calculus is the fractional derivative, 
and some well-known definitions of it are found in the literature, like 
those associated with the names of Riemann-Liouville, Caputo, Weyl, Riesz, etc. 
One common characteristic of these definitions is that they use 
integrals with singular kernels. 
Since the generalization of a concept can sometimes be considered along different 
directions, as in the case of fractional derivative, in recent times new 
definitions of it have been proposed. A recent timeline can be found in 
Chapter I in \cite{eco} and in the references contained therein. 
The fact is that today we seem to have a zoo of definitions of fractional derivative, 
and a classification scheme is certainly welcome.
One such classification was proposed by Teodoro et al. \cite{teodoro}, where 
five different classes were introduced. 

Nevertheless, there are several papers refuting the use of 
the qualification fractional 
for a wide class of those new proposed definitions of fractional derivatives. 
As an example, we mention a recent one like \cite{abdelhakim},  which is a 
sequel of the paper \cite{abdelTenreiro}. In summary, the author 
discusses what he claims to be a flaw in the so-called conformable calculus, 
and he refutes the so-called conformable derivatives 
as proposed, in 2014, by Khalil et al. \cite{khalil} and the whole class of these local derivatives, in the sense they are not a fractional derivative. 
In other words, non-locality is an essential characteristic 
of a fractional derivative \cite{Tarasov}. 

In 2015 there were published two papers; one of them by Caputo-Fabrizio (CF) \cite{CaputoF} proposing a 
new fractional derivative with a \textit{non-singular} kernel, and another one by Losada-Nieto \cite{LosadaN} discussing some 
properties of the so-called CF fractional derivative. 
After these two papers, another definition using a \textit{non-singular} 
kernel was proposed by Atangana-Baleanu (AB) \cite{AtanganaB}, which was used  
by Atangana \cite{Atangana} in the study of the Fisher's reaction-diffusion equation. 
The CF and AB definitions of a fractional
derivative have received considerable interest since their introduction, 
with a combined score of more
than 1000 citations (as of September 2019), being used in the modelling of many 
different problems in terms of
fractional differential equations -- see, for instance \cite{Abro,Tawfik,Kumar, Sun}.


On the other hand, 
recently \cite{capelasvaz} we have reconsidered the use of  
fractional derivatives in the study of the relaxation problem,   
and we have concluded that CF and AB fractional derivatives were
not suitable for the modelling of this problem. The objective
of the present work is to extend that analysis to problems
involving fractional partial differential equations. 
We will show that CF and AB fractional derivatives 
have an intrinsic problem that restricts them from being used to
model problems with initial conditions, which is the 
fact that they use non-singular kernels. Other works have
already identified problems with CF and AB fractional 
derivatives, like, for example, \cite{Stynes,Ortigueira1,Giusti1,Ortigueira2}. 
Our approach can be characterized by its simplicity, as it is based 
on well-known results of Fourier and Laplace transforms.

Possibly the most outstanding example of a fractional partial 
differential equation is the 
version of the diffusion equation with the first order time derivative replaced by a fractional derivative
of order $\alpha$ with $0 < \alpha \leq 1$. Let us see how this equation can 
be obtained. 
In the continuous time random walk (CTRW) approach 
to fractional partial differential equation 
\cite{guide}, we start with the waiting time probability distribution function (PDF) 
$w(t)$ and the jump length PDF $\lambda(x)$. These PDF are related to the PDF $W(x,t)$ of 
being in $x$ at time $t$ through a master equation, which can be solved in the Fourier-Laplace 
space in terms of  the Laplace transform $\tilde{w}(s)$ of $w(t)$ and the Fourier transform $\hat{\lambda}(k)$ 
of $\lambda(x)$.  

In the case of a gaussian jump length PDF with variance $2 \sigma^2$ and a long-tailed waiting time PDF with asymptotic behaviour 
$$
w(t) \sim \frac{\alpha}{\Gamma(1-\alpha)} \frac{\tau^\alpha}{t^{\alpha+1}} , 
\quad (0 < \alpha < 1, t \to \infty) , 
$$
where $\tau$ is a characteristic parameter, we have in the diffusion limit $(s,k) \to (0,0)$ in the Fourier-Laplace space   
that 
$$
\tilde{w}(s) = 1 - \tau^\alpha s^\alpha + \cdots, 
\quad 
\hat{\lambda}(k) = 1 - \frac{\sigma^2}{4} k^2 + \cdots. 
$$
After inversion, from the master equation \cite{guide} it follows the diffusion equation 
\begin{equation}
 \sideset{_{\scriptscriptstyle{\textrm C}}^{}}{_t^{\alpha}}{\operatorname{\mbox D}}[W(x,t)] 
= c^2_\alpha \frac{\partial^2 W(x,t)}{\partial x^2}
\end{equation}
where 
\begin{equation}
\sideset{_{\scriptscriptstyle{\textrm C}}^{}}{_t^{\alpha}}{\operatorname{\mbox D}}[f(t)] = 
G_{1-\alpha}(t)\ast f^\prime(t) = 
\frac{1}{\Gamma(1-\alpha)}\int_0^t \frac{f^\prime(\tau)}{(t-\tau)^\alpha}\, {\textrm{d}}\tau  ,
\end{equation}
with $\ast$ denoting the convolution product, 
is the Caputo fractional derivative. In this expression $G_\nu(t)$ is the Gelfand-Shilov distribution \cite{Dist1,Dist2}, defined as 
\begin{equation}
G_\nu(t) = \begin{cases} {\displaystyle \frac{t^{\nu-1}}{\Gamma(\nu)}H(t)} , & \quad \nu > 0 , \\
{\displaystyle G_{\nu+1}^\prime(t)} , & \quad \nu \leq 0 , 
\end{cases} 
\end{equation}
and satisfying the properties  
\begin{equation}
\mathscr{L}[G_\nu(t)](s) = s^{-\nu} , \qquad 
G_\mu(t)\ast G_\nu(t) = G_{\mu+\nu}(t) , \qquad 
\lim_{\nu\to 0} G_\nu(t) = \delta(t) , 
\end{equation}
where $\mathscr{L}$ denotes the Laplace transform and $H(t)$ the Heaviside step 
function. 

The CF derivative \cite{CaputoF}, denoted by 
$\sideset{_{\scriptscriptstyle{\textrm CF}}^{}}{_t^{\alpha}}{\operatorname{\mbox D}}[f(t)]$, 
and the AB derivative \cite{AtanganaB}, denoted by 
$ \sideset{_{\scriptscriptstyle{\textrm AB}}^{}}{_t^{\alpha}}{\operatorname{\mbox D}}[f(t)]$, are defined  as 
\begin{equation}
\sideset{_{\scriptscriptstyle{\textrm CF/AB}}^{}}{_t^{\alpha}}{\operatorname{\mbox D}}[f(t)]
= \Psi_{{\scriptscriptstyle{\textrm CF/AB}}}(t,\alpha)\ast f^\prime(t) ,
\end{equation}
with 
\begin{equation}
\Psi_{{\scriptscriptstyle{\textrm CF}}}(t,\alpha) = 
\frac{M(\alpha)}{1-\alpha} {\mbox e}^{-\kappa_\alpha(t/\tau) }
\end{equation}
and 
\begin{equation}
\Psi_{{\scriptscriptstyle{\textrm AB}}}(t,\alpha) = 
\frac{M(\alpha)}{1-\alpha} \operatorname{E}_\alpha(-\kappa_\alpha(t/\tau)^\alpha ) , 
\end{equation}
respectively, where $\operatorname{E}_\alpha(\cdot)$ is the Mittag-Leffler function with parameter $\alpha$ \cite{ML,ecomainardivaz}, $0 < \alpha \leq 1$ is the
order of the derivative, $M(\alpha)$ is a normalization, and 
\begin{equation}
\kappa_\alpha = \frac{\alpha}{1-\alpha} . 
\end{equation}
One important feature of these definitions of fractional derivative, as emphasized by their proponents, is that they are non-singular 
for $\alpha \neq 1$, that is, 
their kernel are such that 
\begin{equation}
\label{kernel.limit}
\lim_{t\to 0^+} \Psi_{{\scriptscriptstyle{\textrm CF}}}(t,\alpha) = 
\lim_{t\to 0^+} \Psi_{{\scriptscriptstyle{\textrm AB}}}(t,\alpha) = 
\frac{M(\alpha)}{1-\alpha} . 
\end{equation}
The Caputo kernel, on the other hand, is singular for $0 < \alpha \leq 1$, 
\begin{equation}
\lim_{t\to 0^+} G_{1-\alpha}(t) = +\infty . 
\end{equation}

We will now show that it is precisely the characteristic of having a non-singular kernel that precludes the use of 
these definitions of fractional derivative from being used in 
fractional partial differential equations with initial value problem. In terms
of symbols, let us consider a generic fractional derivative $\sideset{}{_t^{\alpha}}{\operatorname{\mbox D}}[f(t)]$ 
 with a kernel $\Psi(t,\alpha)$, that is, 
\begin{equation}
\label{caputo.type}
\sideset{}{_t^{\alpha}}{\operatorname{\mbox D}}[f(t)] = 
\Psi(t,\alpha)\ast f^\prime(t) . 
\end{equation}
Let us consider the usual diffusion problem of fractional order $\alpha$ \cite{Mainardi} for 
$-\infty < x < \infty$ and $t \geq 0$,  
\begin{equation}
\label{frac.dif.PDE}
\begin{split}
& \sideset{_{\scriptscriptstyle{}}^{}}{_t^{\alpha}}{\operatorname{\mbox D}}[W(x,t)] 
= c_\alpha^2 D^2_x W(x,t) ,  \\[1ex]
& {\displaystyle \lim_{x\to \pm \infty}} W(x,t) = 0 , \\[1ex]
& W(x,0) = \phi(x) .
\end{split}
\end{equation}
Taking the Laplace and Fourier transforms, 
we obtain 
\begin{equation}
\label{eq.aux}
\psi(s,\alpha) s \widehat{\widetilde{W}}(k,s) - \psi(s,\alpha) 
\widetilde{\phi}(k) = - c_\alpha^2 k^2 \widehat{\widetilde{W}}(k,s) ,  
\end{equation} 
where the Laplace and Fourier transforms were denoted by a 
tilde and a hat, respectively, 
that is, $\widetilde{W}(x,s) = \mathscr{L}[W(x,t)]$ and 
$\widehat{W}(k,t) = \mathscr{F}[W(x,t)]$,  and where 
we denoted 
\begin{equation}
\psi(s,\alpha) = \mathscr{L}[\Psi(t,\alpha)](s) , 
\end{equation}
and 
\begin{equation}
\widehat{\phi}(k) = \mathscr{F}[\phi(x)](k) = \widehat{W}(k,0)  . 
\end{equation}
From eq.\eqref{eq.aux} we get  
\begin{equation}
\label{eq.aux.2}
\widehat{\widetilde{W}}(k,s) = \frac{\psi(s,\alpha)}{s\psi(s,\alpha)+ c_\alpha^2 k^2} 
\widehat{\phi}(k) , 
\end{equation}
 
The so-called initial value theorem for Laplace transforms \cite{Laplace} 
says that, if ${\displaystyle \lim_{t\to 0^+}}f(t)$ and $F(s) = \mathscr{L}[f(t)](s)$ exist, 
then 
\begin{equation}
\label{init.value.theorem}
{\displaystyle \lim_{t\to 0^+}}f(t) = \lim_{s\to \infty} s F(s) . 
\end{equation} 
In terms of the notation of our problem, taking $\widehat{W}(k,t)$ as the above 
function $f(t)$, the initial value theorem says that 
\begin{equation}
\lim_{t\to 0^+}\widehat{W}(k,0) = \widehat{\phi}(k) = 
\lim_{s\to \infty} s \widetilde{\widehat{W}}(k,s) . 
\end{equation}
From eq.\eqref{eq.aux.2} we have 
\begin{equation}
\label{eq.aa.19}
\widehat{\phi}(k) = 
\lim_{s\to \infty} s \widetilde{\widehat{W}}(k,s) = 
\lim_{s\to \infty} \frac{s \psi(s,\alpha)}{s\psi(s,\alpha)+ c_\alpha^2 k^2} 
\widehat{\phi}(k) , 
\end{equation}
which  implies, for an arbitrary initial condition, that we must have 
\begin{equation}
\label{eq.aa.20}
\lim_{s\to \infty} \frac{s \psi(s,\alpha)}{s\psi(s,\alpha)+ c_\alpha^2 k^2} = 1. 
\end{equation}
Therefore $\psi(s,\alpha)$ has to satisfy 
\begin{equation}
\label{limit.cond}
\lim_{s\to \infty} [s\psi(s,\alpha)]^{-1} = 0 . 
\end{equation}
The initial value problem can be used again for $\mathscr{L}[\Psi(t,\alpha)] = 
\psi(s,\alpha)$, and the above condition implies therefore that 
\begin{equation}
\lim_{t\to 0^+} \Psi(t,\alpha) = \pm \infty. 
\end{equation}

In conclusion: in order to satisfy an arbitrary initial condition for  
the fractional diffusion equation for $0 < \alpha \leq 1$, 
the kernel $\Psi(t,\alpha)$, as in eq.\eqref{caputo.type}, has to be singular. 
Therefore, the CF and the AB 
definitions are ruled out as candidates for fractional 
derivative models involving an initial value problem. 
For these cases, we have 
\begin{equation}
\psi_{{\scriptscriptstyle{\textrm CF}}}(s,\alpha) = M(\alpha)
\frac{1}{(1-\alpha)s+ \alpha \tau^{-1}} , \qquad 
\psi_{{\scriptscriptstyle{\textrm AB}}}(s,\alpha) = M(\alpha) 
\frac{s^{-1}}{(1-\alpha)+\alpha(s\tau)^{-\alpha}} .
\end{equation}
The condition in eq.\eqref{limit.cond} for these cases is  
\begin{equation}
\lim_{s\to \infty} [s\psi_{{\scriptscriptstyle{\textrm CF}}}(s,\alpha)]^{-1} 
= \frac{1-\alpha}{M(\alpha)}  , \qquad 
\lim_{s\to \infty} [s\psi_{{\scriptscriptstyle{\textrm AB}}}(s,\alpha)]^{-1} 
= \frac{1-\alpha}{M(\alpha)}  , 
\end{equation}
which are satisfied only in the limit $\alpha = 1$, where their kernel are singular, as we see from eq.\eqref{kernel.limit}. 

We would like to remark the use of arbitrary initial conditions in the 
above argument, in particular for obtaining the condition in eq.\eqref{eq.aa.20} 
from eq.\eqref{eq.aa.19}. There is a way to circumvent the eq.\eqref{eq.aa.20} 
to be obtained from eq.\eqref{eq.aa.19}, that is to use as initial 
condition $W(x,0) = \phi(x) = 0$. Then $\widehat{\phi}(k) = 0$ and 
eq.\eqref{eq.aa.19} holds trivially, which implies that 
the limit on the LHS of eq.\eqref{eq.aa.20} can be arbitrary. 
However, the usefulness of problems that can be formulated 
only for vanishing initial conditions does not even deserve 
comments. In addition, a change of variable does not help 
to circumvent the problem when we have arbitrary initial 
condition for a fractional partial differential equation. 
In fact, suppose we make the change of variable $\chi(x,t) = 
W(x,t) - W(x,0) = W(x,t) - \phi(x)$. Then the problem of 
eq.\eqref{frac.dif.PDE} transforms to 
\begin{equation}
\label{frac.dif.PDE.trans}
\begin{split}
& \sideset{_{\scriptscriptstyle{}}^{}}{_t^{\alpha}}{\operatorname{\mbox D}}[\chi(x,t)] 
= c_\alpha^2 D^2_x \chi(x,t) + c_\alpha^2 D^2_x \phi(x),  \\[1ex]
& {\displaystyle \lim_{x\to \pm \infty}} \chi(x,t) = 0 , \\[1ex]
& \chi(x,0) = 0 ,
\end{split}
\end{equation}
that is, a fractional diffusion  equation with a non-homogenous term $c_\alpha^2 D^2_x \phi(x)$. 
Taking as above the Laplace and Fourier transforms, we obtain 
\begin{equation}
\widehat{\widetilde{\chi}}(k,s) = - 
\frac{c_\alpha^2 k^2 \widehat{\phi}(k) s^{-1}}{\phi(s,\alpha) s + k^2 c_\alpha^2} . 
\end{equation}
The initial value theorem eq.\eqref{init.value.theorem} for the Laplace transform gives 
\begin{equation}
\lim_{s\to \infty} \frac{-c_\alpha^2 k^2 \widehat{\phi}(k)}{\psi(s,\alpha) s + k^2 c_\alpha^2} = 
\lim_{t\to 0^+} \widehat{\chi}(k,0) = 0 , 
\end{equation}
where in the last equality we used $\widehat{\chi}(k,0) = \mathscr{L}[\chi(x,0)] = 0$. 
Therefore, the limit on the LHS vanishes only if $\widehat{\phi}(k) = 0$, which 
gives the condition of vanishing initial value we are trying to circumvent, or if 
\begin{equation}
\lim_{s\to \infty} \frac{-c_\alpha^2 k^2}{\psi(s,\alpha) s + k^2 c_\alpha^2} = 0 ,
\end{equation}
which gives as before the condition in eq.\eqref{limit.cond}.

Although we have used the diffusion equation as a prototype 
of a fractional partial differential equation, our analysis
is clearly not dependent on it since it is based on a general property of the Laplace transform. Therefore definitions 
of fractional derivatives based on non-singular kernels cannot  be useful in modelling problems of 
initial value type with fractional differential equations for the very simple reason that they give results
that do not satisfy the initial conditions. Essentially the same conclusion has been 
reached by Stynes in \cite{Stynes} using a different approach. In fact, 
Stynes showed, for a fractional partial differential 
equation of the form 
$$
\sideset{_{\scriptscriptstyle{}}^{}}{_t^{\alpha}}{\operatorname{\mbox D}}[u(x,t)] 
= L_x(x,t)[u(x,t)] + f(x,t) , 
$$
where $L_x(x,t)$ is a second-order differential operator with coefficients
depending on $x$ and $t$, and $f(x,t)$ is a non-homogeneous term, that under certain 
very general conditions, the initial condition $\phi(x) = u(x,0)$ is 
uniquely determined by $L_x(x,0)$, $f(x,t)$ and the boundary condition.  

Finally, it is worthy to mention that recently 
a new fractional derivative closed related to AB one was proposed by Giusti 
and Colombaro \cite{Giusti3,Giusti2} and by 
Zhao and Sun  \cite{Zhao}. It is based on the Prabhakar kernel 
\begin{equation}
\label{kernel.prabhakar}
\Psi_{P}^{\beta,\gamma}(t,\alpha) = \frac{1}{1-\alpha} 
t^{\beta-1} \operatorname{E}_{\alpha,\beta}^\gamma(\lambda t^\alpha) , 
\end{equation}
where $\operatorname{E}_{\alpha,\beta}^\gamma(\cdot)$ is the 
three-parameter Mittag-Leffler function or Prabhakar function \cite{ML,Garra}. 
The Mittag-Leffler function, and therefore the AB fractional derivative, 
corresponds to the particular case $\beta=\gamma=1$. 
The Laplace transform of $\Psi_{P}^{\beta,\gamma}(t,\alpha)$ is 
\begin{equation}
\psi_{P}^{\beta,\gamma}(s,\alpha) = \frac{1}{1-\alpha} 
\frac{s^{\alpha\gamma-\beta}}{(s^\alpha-\lambda)^\gamma} , 
\end{equation}
and therefore 
\begin{equation}
{\displaystyle \lim_{s\to\infty}}[s\psi_{P}^{\beta,\gamma}(s,\alpha)]^{-1} = 
(1-\alpha) {\displaystyle \lim_{s\to\infty}} s^{\beta-1} , 
\end{equation}
and the condition in eq.\eqref{limit.cond} is satisfied for 
\begin{equation}
\beta < 1 . 
\end{equation}
From eq.\eqref{kernel.prabhakar} we see that for $\beta<1$ the 
kernel is singular, that is, 
\begin{equation}
\lim_{t\to 0^+} \Psi_{P}^{\beta,\gamma}(t,\alpha)  = \pm \infty, 
\end{equation}
and therefore models based on it are
expected to give results compatible with the initial conditions 
of the problem, like the case with Caputo derivative.  
In other words, the parameter $\beta$ fixes the 
problem associated with the AB fractional derivative. 
Analytical solution for relaxation models 
based on the Prabhakar derivative have been recently studied \cite{Gorska}.

\medskip

\noindent \textbf{Acknowledgements:} SJ is grateful to CNPq for the financial support.

\end{document}